\documentclass{article} 
\usepackage{amsmath} 
\usepackage{amsfonts}
\usepackage{amsmath}
\usepackage{footnote}
\usepackage{multicol}
\usepackage{booktabs}
\usepackage[flushleft]{threeparttable}
\usepackage[font=small,labelfont=bf]{caption}
\newtheorem{theorem}{Theorem}[section]
\newtheorem{lemma}[theorem]{Lemma}
\newtheorem{proof}[theorem]{Proof}

\newtheorem{definition}[theorem]{Definition}
\newtheorem{example}[theorem]{Example}

\title{Bi-$ g $-frame and  characterizations of bi-$ g $-frame  and
Riesz basis} 
\author{Sayyed Mehrab Ramezani \\ Department of Mathematics, Yasouj University, Yasouj, Iran.\\
m.ramezani@yu.ac.ir} 
\begin{document} 
    \maketitle 
    

    \begin{abstract}
In this paper, we define the concept of the bi-$g$-frame and then show some properties of the bi-$g$-frame. Similarly to bi-$g$-frame, we can define Bi-$g$-Bessel sequences, tight bi-$g$-frame, and the Parseval bi-$g$-frame. Moreover, we define the bi-$g$-frame operator. Finally we present characterizations of bi-$g$-frame and Riesz basis.
\end{abstract}
   \section{Introduction and preliminaries }
Frames for Hilbert space were formally defined by Duffin and Schaeffer \cite{duf} in 1952 while studying some problems in non-harmonic Fourier series.
Recall that for a Hilbert space ${H}$ and a countable index set $J$,
a collection $\lbrace f_{j}\rbrace_{j\in J}\subset {H}$ is called a frame for the Hilbert space
${H}$ if there exist two positive constants $a$,  $b$ such that for all $f\in {H}$
\begin{equation}\label{eq1}
a\Vert f\Vert^2\leq\sum_{j\in J}\vert\langle f , f_{j}\rangle\vert^2\leq b\Vert f\Vert^2;
\end{equation}
$a$ and $b$ are called the lower and upper frame bounds, respectively.
If only the right-hand inequality in \eqref{eq1} is satisfied,  $\lbrace f_{j}\rbrace_{j\in J}$  is called a Bessel sequence for ${H}$ with Bessel bound $b$.
The bounded linear operator $T$ is  defined by
$$T:\ell^{2}(J)\longrightarrow{H},\hspace{1cm}
{T}\lbrace c_{j}\rbrace_{j\in J}=\sum_{j\in J}c_{j}f_{j},$$
which is called the synthesis  operator of $\lbrace f_{j}\rbrace_{j\in J}$.
Moreover, ${T}^{*}f=\lbrace\langle f, f_{j}\rangle\rbrace_{j\in J}$ for all $f\in {H}$. The map $ T^* $ is called the analysis operator of $\lbrace f_{j}\rbrace_{j\in J}$.
 The bounded linear operator $S$
is also defined by
$$S=TT^*: H\longrightarrow H,\hspace{1cm}S(f)=\sum_{j\in J}\langle f, f_{j}\rangle f_j,$$
which is called the frame operator of $\lbrace f_{j}\rbrace_{j\in J}$.
 For more information about
the frames see \cite{chri}.

 Two Bessel sequences $\lbrace f_{j}\rbrace_{j\in J}$ and $\lbrace g_{j}\rbrace_{j\in J}$ are said to be duals for ${H}$ if the following equalities hold
$$f =\sum_{j\in J}\langle f , f_{j}\rangle g_{j}=\sum_{j\in J}\langle f , g_{j}\rangle f_{j},\ \text {for all f}\in {H.} $$
Note that because $S:H\longrightarrow H$ by $S(f)=\sum_{j\in J}\langle f, f_{j}\rangle f_j$ is bijective, self-adjoint and 
$$f=S(S^{-1}f)=\sum_{j\in J}\langle S^{-1}f, f_{j}\rangle f_j=\sum_{j\in J}\langle f, S^{-1}f_{j}\rangle f_j,$$
hence  the sequence $\lbrace S^{-1}f_{j}\rbrace_{j\in J}$ is also a frame by Corollary $1.1.3$ \cite{chri} and it is called the
canonical dual of $\lbrace f_{j}\rbrace_{j\in J}$.
Dual frames are important in reconstructing vectors (or signals) in terms of the
frame elements.
Wenchang Sun \cite{sun}  has provided  characterizations of $ g $-frames and  has proved  that $ g $-frames share many useful
properties with frames. Ramezani and Nazari \cite{mehrab} have gave equivalent conditions for a $ g $-orthonormal basis and characterize all $ g $-Riesz basis for a separable Hilbert, starting with a given $ g $-orthonormal basis. Assuming that 
$\mathcal H$ and $\mathcal K$ are two Hilbert spaces and $\lbrace \mathcal K_{j}\rbrace_{j\in J}$ is a sequence of closed Hilbert subspaces of $\mathcal K$. For each $j\in J,$ $\mathcal{L}\left(\mathcal H,\mathcal K_{j} \right) $ is the collection of all bounded linear operators from $\mathcal H$ to $\mathcal K_{j}$.
A sequence  $\lbrace \Lambda_{j}\in\mathcal{L}(\mathcal H,{\mathcal K}_j),\ {j\in J}\rbrace$ is called  a generalized frame, or simply a
$ g $-frame, for $ \mathcal H $ with respect to $ \lbrace{\mathcal K}_j\rbrace_{j\in J} $ if there are two positive constants $ A $ and $ B $ such that 
\begin{equation}\label{eq2}
A\Vert f\Vert^2\leq\sum_{j\in J}\Vert\Lambda_{j}f\Vert^2\leq B\Vert f\Vert^2,
\hspace{3cm}
(f\in\mathcal H).
\end{equation} 
$A$ and $B$ are called the lower and upper $ g $-frame bounds, respectively.
If only the right-hand inequality in (\ref{eq2}) is satisfied,  $\lbrace \Lambda_{j}\rbrace_{j\in J}$ is called a $ g $-Bessel sequence for $\mathcal{H}$ with respect to $ \lbrace{\mathcal K}_j\rbrace_{j\in J} $ with $ g $-Bessel bound $B$.
The bounded linear operator $T$ is  defined by
$$T:\bigoplus_{j\in J}\mathcal K_{j}\longrightarrow\mathcal{H},\hspace{1cm}
{T}\lbrace g_{j}\rbrace_{j\in J}=\sum_{j\in J}\Lambda_j^* g_{j},$$
which is called the synthesis  operator of $\lbrace \Lambda_{j}\rbrace_{j\in J}$.
Moreover, ${T}^{*}f=\lbrace\langle \Lambda_{j}f\rangle\rbrace_{j\in J}$. The map $ T^* $ is called the analysis operator of $\lbrace \Lambda_{j}\rbrace_{j\in J}$.
 The bounded linear operator $S_\Lambda$
is also defined by
$$S_\Lambda=TT^*:\mathcal H\longrightarrow\mathcal H,\hspace{1cm}S_\Lambda(f)=\sum_{j\in J}\Lambda_j^*\Lambda_jf,$$
which is called the $ g $-frame operator of $\lbrace \Lambda_{j}\rbrace_{j\in J}$.

Weighted and controlled frames have been introduced recently to improve
the numerical efficiency of iterative algorithms for inverting the frame operator
on abstract Hilbert spaces \cite{bal}. By decreasing the ratio of the frame bounds, weighting improves the numerical efficiency of iterative algorithms, such as  the "frame algorithm"\cite{chri} for the inversion of the frame operator.
However, they have been employed earlier in \cite{bog} for spherical
wavelets.
Let $ G L( H)$ be the set of all the bounded operators with a bounded inverse. A frame controlled by the operator $C$ or $C$-controlled frame is a family of
vectors $\lbrace f_{j}\rbrace_{j\in J}\subseteq{H}$, such that there exist two constants $A_c > 0$ and $B_c <\infty$,
satisfying
\begin{equation}\label{eq3}
A_c\Vert f\Vert^2\leq\sum_{j\in J}\langle f , f_{j}\rangle\langle Cf_{j},f\rangle\leq B_c\Vert f\Vert^2;
\end{equation}
for every $f\in{H}$, where $C\in GL({H})$. Every frame is an $I$-controlled frame. Hence
the controlled frames are generalizations of frames. 
The controlled frame operator
$S_c$ is defined by
\begin{equation}
S_c f=\sum_{j\in J}\langle f , f_{j}\rangle Cf_{j}=CS, \hspace{1cm} (f\in{H}),
\end{equation}
where $S$ is the  frame operator of $\lbrace f_{j}\rbrace_{j\in J}$.
The synthesis operator for a $C$-controlled  frame $\lbrace f_{j}\rbrace_{j\in J}$ is defined as follows
\begin{equation*}
T_{c}(\lbrace \alpha_{j}\rbrace_{j\in J})=\sum_{j\in J}\alpha_{j} Cf_{j} =CT,
\end{equation*}
where $T$ is the  synthesis operator of $\lbrace f_{j}\rbrace_{j\in J}$ and
$S_c=T_cT^*$.
 $C$-Controlled  frame $\lbrace f_{j}\rbrace_{j\in J}$ and Bessel sequences$\lbrace g_{j}\rbrace_{j\in J}$ are said to be the $C$-controlled  duals for ${H}$ if the following equality holds.
$$f =\sum_{j\in J}\langle f , g_{j}\rangle Cf_{j},\ \text {for all f}\in {H.} $$
 M. Firouzi Parizi, A. Alijani M. A. Dehghan  \cite{mar} have defined
the concept of biframe that is proposed as a generalization of controlled frames and a special case of pair
frames. A biframe is a pair 
 $ \left( \lbrace f_j \rbrace _{j\in J},\lbrace g_j \rbrace _{j\in J}\right) $
 of sequences in the Hilbert space $ U $, if there
exist positive constants $ c $ and $ d $ such that 
\begin{equation}\label{eq4}
c\Vert f\Vert^2\leq\sum_{j=1}^\infty\langle f,f_j\rangle \langle g_j,f\rangle\leq d\Vert f\Vert^2,
\hspace{3cm}
(f\in U).
\end{equation}

Through the exciting developments
in the biframes and controlled frames,  we introduce the notion of bi-$g$-frames
in Hilbert spaces and show some properties of the bi-$g$-frame then establish a relationship between  bi-$g$-frame and the Riesz basis.

\section{Bi-$ g $-frames and their operator }\label{sec3}
In the following of  this paper, $ \mathcal{U} $ and $ \mathcal{V} $ are two Hilbert spaces and $ \lbrace \mathcal{V}_j:\ j\in J\rbrace $ is a sequence of
subspaces of $ \mathcal{V} $ , where $ J $ is a subset of $ \mathbb{Z} $. $ \mathcal{L}( \mathcal{U},\mathcal{V}_j  ) $  is the collection of all bounded linear operators
from $ \mathcal{U} $  into $ \mathcal{V}_j $.
In this section, we define the concept of bi-$g$-frame  then  show some properties of the bi-$g$-frame.
\begin{definition}\label{def101}
Let $ \lbrace\Lambda_j\in\mathcal{L}( \mathcal{U},\mathcal{V}_j  )\rbrace $ and $ \lbrace\Gamma_j\in\mathcal{L}( \mathcal{U},\mathcal{V}_j  )\rbrace $ are two sequence  for $ \mathcal U $ with respect to $ \mathcal{V}_j $.
A pair $ (\Lambda,\Gamma)=(\lbrace\Lambda_j\rbrace_{j\in J},\lbrace\Gamma_j\rbrace_{j\in J}) $  is called a bi-$ g $-frame for
 $ \mathcal U $ with respect to $ \mathcal{V}_j $ if there exist positive constants $ C $ and $ D $ such that 
\begin{equation}\label{eq5}
C\Vert f\Vert^2\leq\sum_{j\in J}\langle\Lambda_{j}f,\Gamma_{j}f\rangle_{\bigoplus_{j\in J}\mathcal V_{j}}\leq D\Vert f\Vert^2,
\hspace{3cm}
(f\in\mathcal U).
\end{equation} 
we call $C$ and $D$ the lower and upper bi-$g$-frame bounds, respectively. If only the right-hand inequality of (\ref{eq5}) is satisfied, we call
$ (\Lambda,\Gamma) $
 the bi-$g$-Bessel sequence for $\mathcal U$ with respect to
$\lbrace \mathcal V_{j}\rbrace_{j\in J}$ with bi-$g$-Bessel bound $D$. If $C=D=\lambda$, we call $ (\Lambda,\Gamma) $ the tight bi-$g$-frame. Moreover, if $C=D=1$, we call $ (\Lambda,\Gamma) $ the Parseval bi-$g$-frame.
\end{definition}
\begin{example}
Let $\mathcal U$ be a separable Hilbert space and $\lbrace f_{j}\rbrace_{j\in J}$and $\lbrace g_{j}\rbrace_{j\in J}$ be two frames for $\mathcal U$. Let $ \Lambda_{f_j }$ and $ \Gamma_{g_j }$ be
the functional induced by $ f_j $ and $ g_j $ respectively, i.e.,

$$\left\lbrace \begin{array}{l}
\Lambda_{f_j }(f)=\langle f,f_j\rangle \\ 
\\
\Gamma_{g_j }(f)=\langle f,g_j\rangle 
\end{array}\right. $$
A pair $ (\Lambda_f,\Gamma_g)=(\lbrace\Lambda_{f_j}\rbrace_{j\in J},\lbrace\Gamma_{g_j}\rbrace_{j\in J}) $  is  a bi-$ g $-frame for
 $ \mathcal U $ with respect to $ \mathbb{C} $.
\end{example}
Let $ (\Lambda,\Gamma) $ be a bi-$g$-frame for $\mathcal U$ with respect to
$\lbrace \mathcal V_{j}\rbrace_{j\in J}$. We define the bi-$g$-frame operator $ S_{\Lambda,\Gamma} $
as follows: 
\begin{equation}\label{eq6}
S_{\Lambda,\Gamma}:\mathcal U\longrightarrow\mathcal U,\hspace{1cm}S_{\Lambda,\Gamma}(f)=\sum_{j\in J}\Gamma^*_j\Lambda_j f,
\end{equation} 
\begin{theorem}
Let $ (\Lambda,\Gamma) $ be a bi-$g$-frame for $\mathcal U$ with respect to
$\lbrace \mathcal V_{j}\rbrace_{j\in J}$ with bounds $ C $ and $ D $. Then
the following statements are true:
\begin{enumerate}
\item
The operator $ S_{\Lambda,\Gamma} $ is well defined, bounded, positive,  invertible with 
$ \Vert S^{-1}_{\Lambda,\Gamma}\Vert\leq\frac{1}{C} $ and $ S^*_{\Lambda,\Gamma} =S_{\Gamma,\Lambda}  $.
\item
$ (\Lambda,\Gamma) $  is a bi-$g$-frame if and only if $ (\Gamma,\Lambda) $ is a  bi-$g$-frame. 
\end{enumerate}
\begin{proof}
\begin{enumerate}
\item
For each $f\in\mathcal U$ we have
 \begin{align*}
  \vert\langle\sum_{j=n_1}^{n_2}\Gamma^*_j\Lambda_j f,f\rangle\vert&
  = \vert\sum_{j=n_1}^{n_2}\langle\Lambda_j f,\Gamma_jf\rangle\vert\\
  &\leq\sum_{j=n_1}^{n_2}\vert\langle\Lambda_j f,\Gamma_jf\rangle\vert.
\end{align*}
Now we see from \eqref{eq5} that the series in \eqref{eq6} are convergent. Therefore, $ S_{\Lambda,\Gamma} $ is well defined. On the other hands
\begin{align*}
\langle S_{\Lambda,\Gamma}(f),f\rangle&=
  \langle\sum_{j\in J}\Gamma^*_j\Lambda_j f,f\rangle\\
  & = \sum_{j\in J}\langle\Lambda_j f,\Gamma_jf\rangle,
\end{align*}
so
\begin{equation*}
C\Vert f\Vert^2\leq\langle S_{\Lambda,\Gamma}(f),f\rangle\leq D\Vert f\Vert^2,
\hspace{3cm}
(f\in\mathcal U),
\end{equation*} 
and this shows that $ S_{\Lambda,\Gamma} $ is positive  and and bounded operator.
To prove that  $ S_{\Lambda,\Gamma} $ is an invertible operator, we need to show that  $ S_{\Lambda,\Gamma} $ and  $ S^*_{\Lambda,\Gamma} $ are injective and have closed ranges. For this
\begin{align*}
\langle S_{\Lambda,\Gamma}(f),g\rangle&=
  \langle\sum_{j\in J}\Gamma^*_j\Lambda_j f,g\rangle\\
  & = \sum_{j\in J}\langle f,\Lambda^*_j\Gamma_jg\rangle\\
  & = \langle f,\sum_{j\in J}\Lambda^*_j\Gamma_jg\rangle\\
  &=\langle f,S_{\Gamma,\Lambda}(g)\rangle.
\end{align*}
Hence $ S^*_{\Lambda,\Gamma}=S_{\Gamma,\Lambda}  $. By the definition of bi-$g$-frame,  $ S_{\Lambda,\Gamma} $ and  $ S_{\Gamma,\Lambda} $ are injective  and have a closed range as in  the proof of Theorem $ 4.8 $ \cite{mar}. Let $ g\in\mathcal U $ be such that $ \langle S_{\Lambda,\Gamma}(f),g\rangle=0 $
for every $ f\in\mathcal U $. Then we have $ \langle f,S_{\Gamma,\Lambda}(g)\rangle=0 $. This implies that $ S_{\Gamma,\Lambda}(g) =0$ and therefore
$ g=0 $. Hence $  S_{\Lambda,\Gamma}(\mathcal U)=\mathcal U $. Consequently, $ S_{\Lambda,\Gamma} $ is invertible and 
\begin{equation*}
C\Vert f\Vert^2\leq\langle S_{\Lambda,\Gamma}(f),f\rangle\leq\Vert S_{\Lambda,\Gamma}(f)\Vert\Vert f\Vert ,
\hspace{3cm}
(f\in\mathcal U),
\end{equation*}
then
\begin{equation*}
C\Vert f\Vert\leq\Vert S_{\Lambda,\Gamma}(f)\Vert ,
\hspace{3cm}
(f\in\mathcal U),
\end{equation*}
so
\begin{equation*}
C\Vert S^{-1}_{\Lambda,\Gamma}(f)\Vert\leq\Vert f\Vert ,
\hspace{3cm}
(f\in\mathcal U),
\end{equation*}
consequently
\begin{align*}
\Vert S^{-1}_{\Lambda,\Gamma}\Vert\leq\frac{1}{C}.
\end{align*}
\item
Let $ (\Lambda,\Gamma) $ be a bi-$g$-frame for $\mathcal U$ with respect to
$\lbrace \mathcal V_{j}\rbrace_{j\in J}$ with bounds $ C $ and $ D $. Then, for every $ f\in\mathcal U $,
\begin{equation*}
C\Vert f\Vert^2\leq\sum_{j\in J}\langle\Lambda_{j}f,\Gamma_{j}f\rangle_{\bigoplus_{j\in J}\mathcal V_{j}}\leq D\Vert f\Vert^2,
\end{equation*} 
and this means that $ \sum_{j\in J}\langle\Lambda_{j}f,\Gamma_{j}f\rangle_{\bigoplus_{j\in J}\mathcal V_{j}} $ belongs to $ \mathbb{R} $, that is, 
\begin{equation*}
\overline{\sum_{j\in J}\langle\Lambda_{j}f,\Gamma_{j}f\rangle_{\bigoplus_{j\in J}\mathcal V_{j}}}=\sum_{j\in J}\langle\Lambda_{j}f,\Gamma_{j}f\rangle_{\bigoplus_{j\in J}\mathcal V_{j}},
\end{equation*} 
 and therefore 
 \begin{equation*}
{\sum_{j\in J}\langle\Gamma_{j}f,\Lambda_{j}f\rangle_{\bigoplus_{j\in J}\mathcal V_{j}}}=\sum_{j\in J}\langle\Lambda_{j}f,\Gamma_{j}f\rangle_{\bigoplus_{j\in J}\mathcal V_{j}},
\end{equation*} 
 so, $ (\Gamma,\Lambda) $  is a bi-$g$-frame for $\mathcal U$ with respect to
$\lbrace \mathcal V_{j}\rbrace_{j\in J}$ with bounds $ C $ and $ D $.
The converse  statement  is obtained in a similar way.
\end{enumerate}
\end{proof}
\end{theorem}
Below we have a theorem that shows that the reconstruction of elements, which is one of the important achievements of frames, can also be achieved by using bi-$g$-frames.
\begin{theorem}
Let $ (\Lambda,\Gamma) $ be a bi-$g$-frame for $\mathcal U$ with respect to
$\lbrace \mathcal V_{j}\rbrace_{j\in J}$ with bi-$g$-frame operator $ S_{\Lambda,\Gamma} $.
Then, for every $ f\in\mathcal U $, the following reconstruction formulas holds: 
\begin{enumerate}
\item
$ f=\sum_{j\in J}\Gamma^*_j\Lambda_j S^{-1}_{\Lambda,\Gamma}f $
\item
$ f=\sum_{j\in J}(\Gamma_jS^{-1}_{\Gamma\Lambda})^*\Lambda_j f $
\end{enumerate}
\begin{proof}
\begin{enumerate}
\item
$f=S_{\Lambda,\Gamma}S^{-1}_{\Lambda,\Gamma}f=\sum_{j\in J}\Gamma^*_j\Lambda_j S^{-1}_{\Lambda,\Gamma}f $,
\vspace*{1cm}
\item
$ f=S^{-1}_{\Lambda,\Gamma}S_{\Lambda,\Gamma}f=S^{-1}_{\Lambda,\Gamma}\sum_{j\in J}\Gamma^*_j\Lambda_j f=
\sum_{j\in J}S^{-1}_{\Lambda,\Gamma}\Gamma^*_j\Lambda_j f=
\sum_{j\in J}(\Gamma_jS^{-1}_{\Gamma\Lambda})^*\Lambda_j f $.
\end{enumerate}
\end{proof}
\end{theorem}
Now let $ \tilde{\Lambda}_j=\Lambda_j S^{-1}_{\Lambda,\Gamma} $ and $ \tilde{\Gamma}_j=\Gamma_j S^{-1}_{\Gamma,\Lambda} $ . Then the above equalities become
\vspace*{0.5cm}
\begin{enumerate}
\item
$ f=\sum_{j\in J}\Gamma^*_j\tilde{\Lambda}_jf $,
\vspace*{0.5cm}
\item
$ f=\sum_{j\in J}(\tilde{\Gamma}_j)^*\Lambda_j f $.
\end{enumerate}
\vspace*{0.5cm}
We show that $ (\tilde{\Lambda},\tilde{\Gamma}) $ is a bi-$g$-Bessel sequence and furthermore this sequence gives rise to expansion coefficients with the minimal norm.
\begin{lemma}
Let $ (\Lambda,\Gamma) $ be a bi-$g$-frame for $\mathcal U$ with respect to
$\lbrace \mathcal V_{j}\rbrace_{j\in J}$ with bi-$g$-frame operator $ S_{\Lambda,\Gamma} $.
Then, $ (\tilde{\Lambda},\tilde{\Gamma})=(\lbrace\tilde{\Lambda}_j\rbrace_{j\in J},\lbrace\tilde{\Gamma}_j\rbrace_{j\in J}) $ is 
a bi-$g$-Bessel sequence for $\mathcal U$ with respect to
$\lbrace \mathcal V_{j}\rbrace_{j\in J}$ with $g$-Bessel bound $\dfrac{1}{C}$.
\begin{proof}
In fact, for any $ f\in\mathcal U $, we have 
\begin{align*}
\sum_{j\in J}\langle\tilde{\Lambda}_{j}f,\tilde{\Gamma}_{j}f\rangle_{\bigoplus_{j\in J}\mathcal V_{j}}&=\sum_{j\in J}\langle\Lambda_j S^{-1}_{\Lambda,\Gamma}f,\Gamma_j S^{-1}_{\Gamma,\Lambda}f\rangle_{\bigoplus_{j\in J}\mathcal V_{j}}\\
&=\langle\sum_{j\in J}\Gamma^*_j \Lambda_j S^{-1}_{\Lambda,\Gamma}f,S^{-1}_{\Gamma,\Lambda}f\rangle_\mathcal U\\
&=\langle S_{\Lambda,\Gamma} S^{-1}_{\Lambda,\Gamma}f,S^{-1}_{\Gamma,\Lambda}f\rangle_\mathcal U\\
&=\langle f,S^{-1}_{\Gamma,\Lambda}f\rangle_\mathcal U\\
&\leq \dfrac{1}{C}\Vert f\Vert^2.
\end{align*} 
\end{proof}
\end{lemma}

\begin{theorem}
Let $ (\Lambda,\Gamma) $ be a bi-$g$-frame for $\mathcal U$ with respect to
$\lbrace \mathcal V_{j}\rbrace_{j\in J}$ with bi-$g$-frame operator $ S_{\Lambda,\Gamma} $.
Then for any $  g_j\in \mathcal V_{j}$ satisfying $ f=\sum_{j\in J}\Gamma^*_j g_j $ $ \big(f=\sum_{j\in J}{\Lambda}^*_j g_j \big)$ we have 
\begin{align*}
\sum_{j\in J}\Vert g_j\Vert^2=\sum_{j\in J}\langle  g_j ,g_j-\tilde{\Gamma}_{j}f\rangle+
\sum_{j\in J}\langle\tilde{\Lambda}_{j}f,\tilde{\Gamma}_{j}f\rangle.
\end{align*}
\begin{align*}\left( 
\sum_{j\in J}\Vert g_j\Vert^2=\sum_{j\in J}\langle  g_j -\tilde{\Lambda}_{j},g_j\rangle+
\sum_{j\in J}\langle\tilde{\Lambda}_{j}f,\tilde{\Gamma}_{j}f\rangle.\right) 
\end{align*}
\begin{proof}
Suppose $ \lbrace g_j\rbrace_{j\in J}\in \bigoplus_{j\in J}\mathcal V_{j} $ is such that $  f=\sum_{j\in J}\Gamma^*_j g_j $ then
\begin{align*}
\sum_{j\in J}\langle\tilde{\Lambda}_{j}f,\tilde{\Gamma}_{j}f\rangle_{\bigoplus_{j\in J}\mathcal V_{j}}&=
\sum_{j\in J}\langle\tilde{\Lambda}_{j}f,{\Gamma}_{j}S^{-1}_{\Gamma,\Lambda}f\rangle_{\bigoplus_{j\in J}\mathcal V_{j}}\\
&=\langle\sum_{j\in J}{\Gamma}^*_{j}\tilde{\Lambda}_{j}f,S^{-1}_{\Gamma,\Lambda}f\rangle_{\mathcal U}\\
&=\langle f,S^{-1}_{\Gamma,\Lambda}f\rangle_{\mathcal U}\\
&=\langle \sum_{j\in J}{\Gamma}^*_j g_j ,S^{-1}_{\Gamma,\Lambda}f\rangle_{\mathcal U}\\
&=\sum_{j\in J}\langle  g_j ,\tilde{\Gamma}_{j}f\rangle_{\bigoplus_{j\in J}\mathcal V_{j}},
\end{align*}
on the other hands
\begin{align*}
&\sum_{j\in J}\langle  g_j ,g_j-\tilde{\Gamma}_{j}f\rangle+
\sum_{j\in J}\langle\tilde{\Lambda}_{j}f,\tilde{\Gamma}_{j}f\rangle\\
&=\sum_{j\in J}\langle  g_j,g_j\rangle-
\sum_{j\in J}\langle  g_j,\tilde{\Gamma}_{j}f\rangle+
\sum_{j\in J}\langle  g_j ,\tilde{\Gamma}_{j}f\rangle\\
&=\sum_{j\in J}\langle g_j,g_j\rangle\\
&=\sum_{j\in J}\Vert g_j\Vert^2,
\end{align*}
therefore
\begin{align*}
\sum_{j\in J}\Vert g_j\Vert^2=\sum_{j\in J}\langle  g_j ,g_j-\tilde{\Gamma}_{j}f\rangle+
\sum_{j\in J}\langle\tilde{\Lambda}_{j}f,\tilde{\Gamma}_{j}f\rangle.
\end{align*}
The second part is proved in the same way.
\end{proof}
\end{theorem}
\section{Characterizations of bi-$ g $-frame  and
Riesz basis }
We get characterizations of 
bi-$g$-Bessel sequence, tight bi-$g$-frame  and bi-$ g $-frame.
\begin{theorem} \label{th4}
Let $ {\Lambda}_{j},\ {\Gamma}_{j}\in \mathcal{L}(\mathcal U,\mathcal V_{j}) $ and $ u_{j,k}= \Lambda^*_je_{j,k}$ , $ v_{j,k}= \Gamma^*_je_{j,k}$ where in, $ \lbrace e_{j,k},\ k\in\mathbb{K}_j\rbrace $ is an orthonormal basis for $\mathcal V_{j}  $and $ \mathbb{K}_j $ is a subset of $ \mathbb{Z} $. Then $ (\Lambda,\Gamma) $ 
 is a bi-$g$-frame (respectively bi-$g$-Bessel sequence, tight bi-$g$-frame) for $\mathcal U$ with respect to
$\lbrace \mathcal V_{j}\rbrace_{j\in J}$ if and only if $ \left( \lbrace u_{j,k} \rbrace _{ {j\in J}, k\in\mathbb{K}_j},\lbrace v_{j,k} \rbrace _{ {j\in J}, k\in\mathbb{K}_j}\right) $ is a biframe (respectively
biBessel sequence, tight biframe) for $ \mathcal U $. 
\begin{proof}
For any  $ f\in\mathcal U $ we have
\begin{align*}
\sum_{j\in J}\langle\Lambda_{j}f,\Gamma_{j}f\rangle_{\bigoplus_{j\in J}\mathcal V_{j}}&=
\sum_{j\in J}\left\langle \sum_{k\in\mathbb{K}_j}\langle \Lambda_{j}f, e_{j,k}\rangle e_{j,k},
\sum_{l\in\mathbb{K}_j}\langle \Gamma_{j}f, e_{j,l}\rangle e_{j,l}\right\rangle _{\bigoplus_{j\in J}\mathcal V_{j}}\\
&=\sum_{j\in J}\sum_{k\in\mathbb{K}_j}\sum_{l\in\mathbb{K}_j}\left\langle \langle f, \Lambda^*_{j}e_{j,k}\rangle e_{j,k},
\langle f, \Gamma^*_{j}e_{j,l}\rangle e_{j,l}\right\rangle _{\mathcal U}\\
&=\sum_{j\in J}\sum_{k\in\mathbb{K}_j}\sum_{l\in\mathbb{K}_j}\left\langle \langle f, u_{j,k}\rangle e_{j,k},
\langle f, v_{j,l}\rangle e_{j,l}\right\rangle _{\mathcal U}\\
&=\sum_{j\in J}\sum_{k\in\mathbb{K}_j} \langle f, u_{j,k}\rangle \cdot
\langle v_{j,k},f\rangle.
\end{align*}
Hence $ (\Lambda,\Gamma) $  is a bi-$g$-frame (respectively bi-$g$-Bessel sequence, tight bi-$g$-frame) for $\mathcal U$ with respect to
$\lbrace \mathcal V_{j}\rbrace_{j\in J}$ if and only if $ \left( \lbrace u_{j,k} \rbrace _{ {j\in J}, k\in\mathbb{K}_j},\lbrace v_{j,k} \rbrace _{ {j\in J}, k\in\mathbb{K}_j}\right) $ is a biframe (respectively
biBessel sequence, tight biframe) for $ \mathcal U $.
\end{proof}
\end{theorem}
The next theorem expresses the dependence of the $ g $-Riesz basis, which are in the form of a bi-$g$-frame.
\begin{theorem} \label{th5}
Let $ (\Lambda,\Gamma) $ be a bi-$g$-frame for $\mathcal U$ with respect to
$\lbrace \mathcal V_{j}\rbrace_{j\in J}$ with bi-$g$-frame operator $ S_{\Lambda,\Gamma} $.
Then $ \Lambda=\lbrace\Lambda_j\rbrace_{j\in J}$ is a $ g $-Riesz basis
for $\mathcal U$ with respect to
$\lbrace \mathcal V_{j}\rbrace_{j\in J}$ if and only if $ \Gamma=\lbrace\Gamma_j\rbrace_{j\in J}$ is a $ g $-Riesz basis for for $\mathcal U$ with respect to
$\lbrace \mathcal V_{j}\rbrace_{j\in J}$. 
\begin{proof}
Let $ (\Lambda,\Gamma) $ be a bi-$g$-frame for $\mathcal U$ with respect to
$\lbrace \mathcal V_{j}\rbrace_{j\in J}$ then by Theorem \ref{th4} $ \left( \lbrace u_{j,k} \rbrace _{ {j\in J}, k\in\mathbb{K}_j},\lbrace v_{j,k} \rbrace _{ {j\in J}, k\in\mathbb{K}_j}\right) $ is a biframe where in $ \lbrace e_{j,k},\ k\in\mathbb{K}_j\rbrace $ is an orthonormal basis for $\mathcal V_{j}  $, $ \mathbb{K}_j $ is a subset of $ \mathbb{Z} $,  $ u_{j,k}= \Lambda^*_je_{j,k}$ and $ v_{j,k}= \Gamma^*_je_{j,k}$.
Now if $ \lbrace\Lambda_j\rbrace_{j\in J}$ is a  $ g $-Riesz basis, by Theorem 3.1 \cite{sun},
the sequence $ \lbrace u_{j,k} \rbrace _{ {j\in J}, k\in\mathbb{K}_j} $ is a Riesz basis. Since  $ \left( \lbrace u_{j,k} \rbrace _{ {j\in J}, k\in\mathbb{K}_j},\lbrace v_{j,k} \rbrace _{ {j\in J}, k\in\mathbb{K}_j}\right) $ be a biframe,   by Theorem 4.8 \cite{mar}, the sequence $ \lbrace v_{j,k} \rbrace _{ {j\in J}, k\in\mathbb{K}_j} $ is also a Riesz basis. Once again, from Theorem 3.1 \cite{sun}, we conclude that  $ \Gamma$ is a $ g $-Riesz basis for for $\mathcal U$ with respect to
$\lbrace \mathcal V_{j}\rbrace_{j\in J}$. 
It can be proved in the same way as the converse of the theorem. 
\end{proof}
\end{theorem}


\end{document}